\documentclass[a4paper,12pt]{article}
\usepackage{theorem,amsmath,amssymb}
\usepackage{epsfig,psfrag}
\usepackage{multicol}

\newcommand{\mysection}{\setcounter{equation}{0} \section}

\newtheorem{theo}{Theorem}[section]
\newtheorem{prop}[theo]{Proposition}
\newtheorem{lemma}[theo]{Lemma}

\newtheorem{fact}[theo]{Fact}
\newtheorem{deff}[theo]{Definition}

\def\a{\alpha}
\def\b{\beta}

\def\d{\delta}
\def\e{\varepsilon}

\def\E{\mathbb{E}}
\def\o{\omega}

\def\l{\lambda}
\def\N{\mathbb{N}}
\def\P{\mathbb{P}}

\def\R{\mathbb{R}}

\def\Z{\mathbb{Z}}

\def\tb{\tau_{-1}}
\def\c{c_{1,n}}
\def\cc{c_{3,n}}

\def\cccc{c_{2,n}}
\def\dd{c_{4,n}}
\def\ddd{c_{5,n}}
\def\dddd{\cccc}
\def\ee{c_{6,n}}
\def\eee{\ddd}

\def\V{\widetilde{V}}
\def\VV{\widehat{V}}

\def\EA{E_{n}}
\def\EE{E_{1,n}}
\def\EEE{E_{2,n}}
\def\Ff{E_{3,n}}
\def\FF{E_{4,n}}
\def\FFF{E_{5,n}}
\def\Gg{E_{7,n}}
\def\GG{E_{2,n}'}
\def\GGG{E_{3,n}'}
\def\Hh{E_{6,n}'}

\def\po{P_{\omega}}

\def\eo{E_{\omega}}
\def\Iq{I^q_{\eta}}

\def\ra{\rightarrow}

\def\v1{v_{1}}

\def\j1{J_{\Lambda}^<}

\def\pas{$\mathbb{P}$--a.s. }

\def\k{\kappa}
\def\jk{J_{\k}}

\def\gk{\Gamma_{\k}}
\def\ik{I_{\k}}
\def\fv{\phi_{v_{\k}}}

\title{Some properties of the rate function of quenched large
deviations for random walk in random environment}

\author{Alexis Devulder
\footnote{Laboratoire de Probabilit\'es et Mod\`eles Al\'eatoires,
Universit\'e Paris VI, $4$ Place Jussieu, Case 188, F-$75252$
Paris Cedex $05$, France. E-mail: devulder@ccr.jussieu.fr.}}

\begin{document}


\maketitle

\begin{abstract}
In this paper, we are interested in some questions of Greven and
den Hollander \cite{GD2} about the rate function $\Iq$ of quenched
large deviations for random walk in random environment. By
studying the hitting times of RWRE, we prove that in the recurrent
case, $\lim_{\theta\to 0^+}(\Iq)''(\theta )=+\infty$, which gives
an affirmative answer to a conjecture of Greven and den Hollander
\cite{GD2}. We also establish a comparison result between the rate
function of quenched large deviations for a diffusion in a drifted
Brownian potential, and the rate function for a drifted Brownian
motion with the same speed.
\end{abstract}


\bigskip
\noindent\textbf{AMS $(2000)$ Classification:} 60K37, 60F10,
60J60.


\mysection{Introduction}\label{intro}
\subsection{Presentation of the model}
We consider a collection of independent and identically
distributed random variables $(\o_i)_{i\in\Z}$. A realization of
these variables is called an environment. Given an environment
$\o:=(\o_i)_{i\in\Z}$, we consider the random walk
$(S_n)_{n\in\N}$ defined by $S_0=0$ and
\begin{equation*}
\po(S_{n+1}=k|S_n=i) =\left\{ {\begin{array}{ll}
 \o_i  & \text{ if } k=i+1,\\
 1-\o_i & \text{ if } k=i-1,\\
 0 & \text{ otherwise.}
\end{array}}
\right.
\end{equation*}
The process $(S_n)_{n\in\N}$ is called a random walk in random
environment, abbreviated RWRE. This model has many applications in
physics, see for example Hughes \cite{Hug}. Let $\eta$ denote the
law of $(\o_i)_{i\in\Z}$. We call $\po$ the quenched law, whereas
$\P(.):=\int \po(.)\eta(d\o)$ is the annealed law. For technical
reasons, we assume that there exists an $\e_0>0$ such that
\begin{equation}
\eta(\e_0\leq\o_0\leq 1-\e_0)=1.\label{lip}
\end{equation}
For $i\in\Z$, let $\rho_i=\frac{1-\o_i}{\o_i}$. Solomon \cite{S1}
proved that the RWRE $(S_n)_{n\in\N}$ is \pas recurrent if and
only if
\begin{equation}
\int(\log\rho_0)\eta(d\o)=0.\label{rec}
\end{equation}
In order to avoid the degenerate case of simple random walk, we
assume in the following that
\begin{equation}
\text{Var}(\log \rho_0):=\sigma^2>0.\label{var}
\end{equation}
Sinai \cite{S3} showed that in the recurrent case, the random
environment considerably slows down the walk. More precisely, he
proved that if \eqref{rec} and \eqref{var} are satisfied, there
exists a nondegenerate non--Gaussian random variable $b_{\infty}$
such that
\begin{equation}
\sigma^2 \frac{S_n}{(\log
n)^2}\overset{\mathcal{L}}{\underset{n\to+\infty}{\longrightarrow
}}b_{\infty},\label{sinai}
\end{equation}
where $\overset{\mathcal{L}}{\longrightarrow}$ denotes convergence
in law under $\P$.

It is moreover known (Solomon~\cite{S1}) that the RWRE
$(S_n)_{n\in\N}$ satisfies a law of large numbers: there exists
$v\in ]-1,1[$ such that $\lim_{n\to\infty} S_n/n=v$ $\P$--a.s. In
addition, $v$ is strictly positive if and only if $\int
\rho_0\eta(d\o)<1$.

The RWRE $(S_n)_{n\in\N}$ satisfies furthermore a quenched large
deviation principle with deterministic convex rate function $\Iq$
(see Greven and den Hollander~\cite{GD2}). This means there exists
a nonnegative convex function $\Iq$ such that $\eta$--a.s. for any
measurable set $A$,
\begin{eqnarray*}
    \liminf_{n\to\infty}\frac{1}{n}\log\po\left(\frac{S_n}{n}\in A\right)
 &\ge& -\inf_{x\in A^{\circ}}\Iq(x),
    \\
    \limsup_{n\to\infty}\frac{1}{n}\log\po\left(\frac{S_n}{n}\in A\right)
 &\le& -\inf_{x\in \overline{A}}\Iq(x),
\end{eqnarray*}
where $A^{\circ}$ denotes the interior of $A$ and $\overline{A}$
is the closure of $A$.

For more details on RWRE, we refer to Zeitouni \cite{Z1}.


\subsection{Results}
In this paper, we are interested in some questions raised by
Greven and den Hollander about quenched large deviations for RWRE.
First, we answer their Open problem 2 (see \cite{GD2}, p. 1389;
see also den Hollander \cite{H3} p. 80), and prove that

\begin{theo}\label{th1}
Under \eqref{lip}, \eqref{rec} and \eqref{var}, the rate function
$\Iq$ for quenched large deviations of the RWRE satisfies
\begin{equation}
\lim_{\theta\ra 0^+} (\Iq)''(\theta)=+\infty.\label{pb1}
\end{equation}
\end{theo}

This is what Greven and den Hollander conjectured. Observe that
this result is coherent with the subdiffusive behaviour of Sinai's
walk \eqref{sinai}. We mention that the corresponding problem for
Brox--type diffusions (see Brox~\cite{B3}), for which the rate
functions can be explicitly computed, has already been solved by
Taleb (see \cite{Ta}).


In order to prove Theorem \ref{th1}, it is useful to study the
hitting times of $(S_n)_{n\in\N}$. Let us define, for $a\in\Z$,
\begin{equation*}
\tau_a:=\inf\{n>0, S_n=a\}.
\end{equation*}
We show the following estimate:

\begin{prop}\label{moments}
For each $\alpha\in \R_+^*$,
\begin{equation*}
\E(\tau_1^\a
e^{-r\tau_1})=\left(\frac{1}{r}\right)^{\a+o(1)},\qquad r\to 0^+.
\end{equation*}
\end{prop}

\bigskip
We are also interested in Open problem 3 of Greven and den
Hollander (\cite{GD2}, p. 1389): they conjectured that in the case
$\int \rho_0 \eta(d\o)<1$ (i.e., $v>0$), the quenched rate
function $\Iq$ of the RWRE is smaller than
the rate function of the simple random walk on $\Z$ with the same
speed $v$. That is, they conjectured that $\forall x>v,\
\Iq(x)<\widehat{I}_{\langle \rho\rangle }(x)$, where
$\widehat{I}_{\langle\rho\rangle}$ is the rate function of a usual
nearest neighbour random walk with speed $v$.

Unfortunately, we have not been able to answer this question, but
we solve the corresponding problem for Brox--type diffusion (see
Brox, \cite{B3}). For $\k\geq 0$, we define the random potential
\begin{equation*}
W_{\k}(x):=W(x)-\frac{\k}{2}x,
\end{equation*}
where $(W(x),x\in\R)$ is a standard two-sided Brownian motion. We
consider a diffusion~$X$ in the random potential $W_{\k }$, which
is defined as the solution to the formal stochastic differential
equation
\begin{equation*}
\left\{
\begin{array}{l}
\textnormal{d}X(t)=\textnormal{d}\b(t)-\frac{1}{2}W'(X(t))\textnormal{d}t,\\
X(0)=0,
\end{array}\right.
\end{equation*}
where $(\b(t), t\geq 0)$ is a Brownian motion independent of $W$.
More precisely, $X$ is a diffusion process whose conditional
generator given $W_{\k}$ is
\begin{equation*}
\frac{1}{2}e^{W_{\k}(x)}\frac{\partial}{\partial
x} \left(e^{-W_{\k}(x)}\frac{\partial}{\partial x}\right).
\end{equation*}
This diffusion can be considered as the continuous time analogue
of RWRE and share many properties with it. See for example Shi
\cite{S2} for the relations between these two processes. For
instance, Kawazu and Tanaka \cite{KT1} established a law of large
numbers for $X$. That is, $\lim_{t\to\infty} X_t/t=v_{\k}$, where
$v_{\k}=\frac{(\k-1)^+}{4}$ is $>0$ if and only if $\k>1$.
Moreover, Taleb \cite{Ta} proved that $X$ satisfies quenched and
annealed large deviation principles. Let $\jk$ denote the rate
function of quenched large deviations of~$X$ (see \eqref{quotient}
below for more details). We compare $\jk$ with the function
\begin{equation*}
J_{v_{\k}}^B(x):=\frac{1}{2}(x-v_{\k})^2,
\end{equation*}
which is the rate function of large deviations of the drifted
Brownian motion $(B_t+v_{\k} t,\ t\in\R_+)$. We prove

\begin{theo}\label{comptaux}
If $\k>1$, then
\begin{equation*}
\forall x>v_{\k},\qquad \jk(x)< J_{v_{\k}}^B(x).
\end{equation*}
\end{theo}

\bigskip
Interestingly, we obtain as a by-product an inequality for the
modified Bessel functions
which might be new:

\begin{prop}\label{prop14}
Let $K_{\nu}$ be the modified Bessel function of index $\nu$. We
have,

\begin{equation*}
\forall \nu>0,\ \forall y>0,\qquad
\frac{K_{\nu}(y)}{K_{\nu+1}(y)}<\frac{1}{y}\left(\sqrt{y^2
+\nu^2}-\nu\right).
\end{equation*}
\end{prop}

\bigskip
The rest of the paper is organized as follows: in Section
\ref{milieu} we build environments $\EA$ for which the hitting
time of~$-1$ by $(S_n)_{n\in\N}$, denoted by $\tau_{-1}$, will be
large. We give an estimation of $\tau_{-1}$ for $\o\in \EA$ in
Section~\ref{proba}. In Section \ref{simpl}, we prove Theorem
\ref{th1} and Proposition \ref{moments}. Finally,
Section~\ref{sectbrox} is devoted to the proofs of
Theorem~\ref{comptaux} and Proposition~\ref{prop14}.


\mysection{Construction of the event {$\hbox{\slshape\bfseries
E}_{\hbox{\slshape\bfseries \small n}}$} } \label{milieu}

In this section we build a set of environments $\EA$, such that
$\P(\EA)$ is not ``too small'' and that for $\o\in \EA$,
$\tau_{-1}$ is almost $n$ (we prove this last assertion in Section
\ref{proba}).

Throughout Sections \ref{milieu} and \ref{proba}, we fix an
$\e>0$. The constants $C_i$, $0\leq i\leq 10,$ depend only on
$\eta$ and $\e$, whereas $\theta$ and $\delta$ depend only on
$\eta$. The events $\EA$, $E_{i,n}$ and $E_{i,n}'$ depend on
$\eta$ and $\e$, but we omit to write $\e$.

We give some notation in Subsection \ref{sub20}. Subsection
\ref{sub21} is devoted to the construction of $\EA$. We give an
estimation of $\P(\EA)$ in Subsection \ref{sub22}, and study some
of the properties of $\o\in \EA$ in Subsection~\ref{sub23}.


\subsection{Some notation}\label{sub20}

We define the potential $V$ as follows:

\begin{deff}Let
\begin{equation*} V(n):=\sum_{i=1}^n\log\rho_i=\sum_{i=1}^n \log\frac{1-\o_i}{\o_i},
\qquad n\in\Z,
\end{equation*}
where by convention, $\sum_{i=1}^0 x_i=0$ and $\sum_{i=1}^n x_i
=-x_0-x_{-1}-\dots -x_{n+1}$ if $n$ is (strictly) negative.
\end{deff}

We define a valley for the potential (see Sinai, \cite{S3}):
\begin{deff}\label{defval}
Let $a<m<b$. $(a,m,b)$ is a valley if
\begin{equation*}
\begin{array}{ll}
\forall a\leq i\leq m, & \quad V(m)\leq V(i)\leq V(a), \nonumber\\
\forall m\leq i \leq b, & \quad V(m)\leq V(i)\leq V(b).\nonumber
\end{array}
\end{equation*}
Its depth is defined as $\min\{V(a)-V(m),V(b)-V(m)\}$.
\end{deff}


\subsection{Building {$\hbox{\slshape\bfseries
E}_{\hbox{\slshape\bfseries \small n}}$}}\label{sub21}

In this subsection, we build a valley $(0, m_n, b_n)$ for the
potential $V$, so that the RWRE will stay for a ``good'' amount of
time in this valley with ``large probability''.

As $\int(\log\rho_0)\eta(d\o)=0$ and $\sigma>0$, there exists a
real number $\d>0$ such that
\begin{equation*}
\P(-2\delta\leq \log \rho_0\leq -\delta):=\exp(-\theta)>0.
\end{equation*}
Now we set
\begin{equation*}
\begin{array}{rclrcl} \e'& := &\e\d,
& \c & := & \lfloor \e \log n\rfloor,\\

 \cccc & := &  \lfloor \log n\rfloor^2,
& \cc & := & \d \c, \\

\dd & := &(1-10\e') \log n,&
\beta & := & \frac{1-9\e'}{1-10\e'},\\

\ddd & := & \frac{\e'}{2} \log n, &
\ee & := & 2 \log n.\\
\end{array}
\end{equation*}

\begin{figure}[htbp]
\begin{center}
\setlength{\unitlength}{1cm}
\begin{picture}(12.5,7)
\footnotesize \psfrag{y}{\footnotesize $j$}
\psfrag{x}{\footnotesize $i$} \psfrag{c1+c2}{\footnotesize
$\c+\cccc$} \psfrag{c1+2c2}{\footnotesize $\c+2\cccc$}
\psfrag{c1}{\footnotesize $\c$} \psfrag{c5}{\footnotesize $\ddd$}
\psfrag{-1}{\footnotesize $-1$}\psfrag{-c3}{\footnotesize
$-\cc$}\psfrag{-2c3}{\footnotesize $-2\cc$}
\psfrag{-c4}{\footnotesize $V(\c)-\dd$}
\psfrag{-bc4}{\footnotesize $V(\c)-\beta\dd$}
\psfrag{c6}{\footnotesize $V(\c+\cccc)+\ee$}
\psfrag{2c6}{\footnotesize $V(\c+\cccc)+2\ee$}
\psfrag{m}{\footnotesize $m_n$} \psfrag{b}{\footnotesize $b_n$}

\put(0,0){\includegraphics[width=12.5cm,height=7cm]{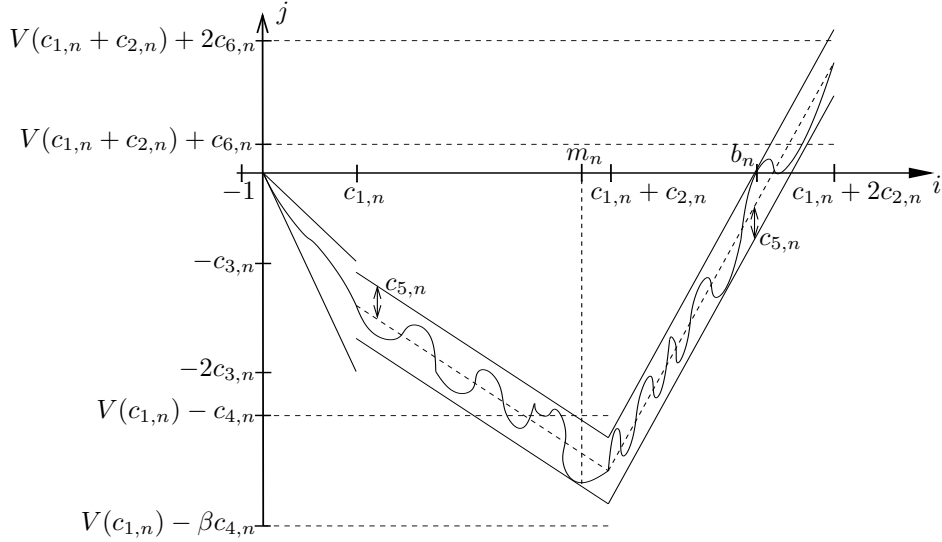}}
\end{picture}
\end{center}\caption{Example of the potential $V$ for a ``good'' environment $\o\in E_n$}
\label{figure1}\end{figure}

\noindent For $i\in\Z$, define $\V(i)=V(i+\c)-V(\c)$ and
$\VV(i)=\V(i+\cccc)-\V(\cccc)$. We consider

\begin{eqnarray*}
\EE & := & \left\{\forall 0\leq i\leq \c, \qquad -2\d i\leq
V(i)\leq -\d i \right\},\\
\EEE & := & \{\V(\cccc)\in [-\b\dd,-\dd]\},\\
\Ff & := & \left\{\forall 0< i \leq \cccc,\qquad
\left|\V(i)-\frac{i}{\cccc}
\V(\cccc)\right|\leq \ddd \right\},\\
\FF & := & \{\VV(\dddd)\in[\ee,2\ee]\},\\
\FFF & := & \left\{\forall 0< i\leq \dddd,\qquad
\left|\VV(i)-\frac{i}{\dddd}\VV(\dddd)\right|\leq \eee \right\}.
\end{eqnarray*}
Finally, let
\begin{equation*}
\EA:=\EE\cap \EEE\cap \Ff\cap \FF\cap \FFF.
\end{equation*}
When $\o\in \EA$, we say the environment is ``good''. On $\EE$,
the potential $V(i)$ decreases almost linearly for $i\in[0,\c]$
(this will enable the walk $(S_n)_{n\in\N}$ to go quickly to $\c$
before hitting $-1$ with large probability). On $\EEE\cap \Ff$,
$V$ stays within a tunnel of height $2\ddd$, and sinks to
$V(\c)-\dd$. On $\FF\cap \FFF$, $V$ stays within another tunnel of
height $2\ddd$, and moves up to positive values. These comments on
$E_n$ are represented in the Figure \ref{figure1} ($b_n$ and $m_n$
are defined in Subsection~\ref{sub23}).

\bigskip

\subsection{Probability of {$\hbox{\slshape\bfseries
E}_{\hbox{\slshape\bfseries \small n}}$}}\label{sub22}
\begin{lemma}\label{en}
There exists a constant $C_0>0$ which depends only on $\eta$ and
$\e$, such that for $n$ large enough,
\begin{equation*}
P(\EA)\geq C_0 n^{-\theta\e}.
\end{equation*}
\end{lemma}

\bigskip
\noindent {\bf Proof:} First, observe that
\begin{eqnarray}
\P(\EE) & \geq & \P\left(\forall 1\leq i\leq \c,\quad
-2\d\leq\log\rho_i\leq-\d \right)\nonumber\\
& \geq & \exp(-\theta {\e\log n})\geq n^{-\theta\e}.\label{E2}
\end{eqnarray}

According to the Koml\'os--Major--Tusn\'ady strong approximation
theorem (see \cite{KMT}), possibly in an enlarged probability
space, there exists a coupling for  $\o$ and a standard Brownian
motion $W$, and (strictly) positive constants $C_1$, $C_2$ and
$C_3$ such that for all $N\geq 1$,
\begin{equation*}
\P\left(\sup_{1\leq i\leq N} |\V(i)-\sigma W(i)|\geq C_1\log
N\right) \leq \frac{C_2}{N^{C_3}}.
\end{equation*}
Define
\begin{equation*}
\Hh:=\left\{\sup_{1\leq i\leq 3(\log n)^2} |\V(i)-\sigma W(i)|\geq
C_1\log [3(\log n)^2]\right\}.
\end{equation*}
We have,
\begin{equation*}
\P(\Hh)\leq \frac{C_2}{(3(\log n)^2)^{C_3}}.
\end{equation*}
We then consider the following events:
\begin{eqnarray*}
\GG & := & \left\{\sigma W(\cccc)\in\left[-\b\dd+\frac{\e'}{4}\log
n,
-\dd-\frac{\e'}{4} \log n \right]\right\},\\
\GGG & := & \left\{\forall 0\leq t\leq \cccc,\qquad \left|\sigma
W(t)-\frac{t}{\cccc}\sigma W(\cccc)\right|\leq
\ddd-\frac{\e'}{4}\log n\right\}.
\end{eqnarray*}
We notice by scaling that there exists $C_4>0$ such that
$\P[\GG\cap \GGG]\geq 2 C_4$ for $n$ large enough. Since
$\log[3(\log n)^2]=o(\log n)$, we have for large $n$,
\begin{eqnarray}
\P(\EEE\cap \Ff) & \geq & \P[\GG\cap \GGG\cap
(\Hh)^c] \nonumber \\
& \geq  & \P(\GG\cap \GGG)-\frac{C_2}{[3(\log n)^2]^{C_3}}\nonumber \\
& \geq & C_4.\label{c10}
\end{eqnarray}
Similarly, there exists a constant $C_5>0$ such that
\begin{equation}\label{c11}
\P(\FF\cap \FFF)\geq C_5
\end{equation}
for $n$ large enough. Since  $\EE$, $\EEE\cap \Ff$ and $\FF\cap
\FFF$ are independent, we obtain Lemma \ref{en} by combining
\eqref{E2}, \eqref{c10} and \eqref{c11}. \hfill$\Box$


\subsection{Properties of a ``good'' environment}\label{sub23}
Let $\o\in E_n$. We define the integers $b_n$ and $m_n$ such that
\begin{eqnarray*}
b_n & := &\inf\{k\in\N, \ k>0, \quad V(k) \geq 0\}, \\
m_n & := & \inf\{k>0 ,\quad V(k)=\inf_{0\leq\ell\leq b_n}
V(\ell)\}.
\end{eqnarray*}
Note  that $(0,m_n,b_n)$ is a valley  (in the sense of Definition
\ref{defval}) with depth $-V(m_n)$, and that
\begin{eqnarray*}
V(m_n) & \in & [-2\cc-\b\dd-\ddd, -\cc-\dd],\\
m_n & \in & \left[\c+\cccc-\frac{\ddd\cccc}{\dd},
\c+\cccc+\frac{\eee\dddd}{\ee}\right].
\end{eqnarray*}
In particular, we have for $\e'$ small enough and $n$ large
enough, \begin{equation}\label{mn}
\begin{array}{rcccl}
-\d+(1-9\e') \log n   & \leq & -V(m_n) & \leq &  (1-6\e') \log n, \\
(1-\e') (\log n)^2 & \leq & m_n      & \leq & (1+\e') (\log n)^2 .
\end{array}
\end{equation}


\mysection{Probability that
$\tau_{-1}$ has a
``good'' length }\label{proba}

This section is devoted to the proof of the following result:

\begin{lemma}\label{encadrement}
There exists a constant $C_6>0$, depending only on $\eta$ and
$\e$, such that for all large $n$,
\begin{equation*}
\forall\o\in \EA,\qquad
\po\left(n^{1-10\e'}\leq\tau_{-1}<n\right)\geq C_6.
\end{equation*}
\end{lemma}

In Subsection \ref{sub31}, we show that when $\o\in \EA$, with a
large quenched probability, the RWRE goes quickly to the bottom
$m_n$ of the valley $(0,m_n,b_n)$ without hitting $-1$. In
Subsection \ref{sub32}, we prove that with a large quenched
probability, after hitting $m_n$, the RWRE stays in $\N$ during
almost $n$ units of time and then hits $-1$ for the first time.


\subsection{Going to the bottom
{$\hbox{\slshape\bfseries m}_{\hbox{\slshape\bfseries \small n}}$}
of the valley}\label{sub31}

\begin{lemma} \label{lemma32}There exists a constant $C_7>0$, depending
only on $\eta$ and $\e$, such that
$$\forall \o\in \EA,\qquad \po(
\tau_{m_n}<\tau_{-1} )\geq C_7.$$
\end{lemma}

\bigskip
\noindent {\bf Proof:} Let $\o\in \EA$. Since $\EE\subset \EA$,
\begin{equation}\label{minorationA1}
\sum_{i=0}^{\c-1}\exp(V(i))\leq \sum_{i=0}^{+\infty}\exp(-\d i)
\leq \frac{1}{1-e^{-\d}}.
\end{equation}
Furthermore,
\begin{equation*}
\forall \c\leq i\leq m_n,\qquad V(i) \leq  -\d\c+\ddd\leq
-\frac{\e'}{2}\log n+\d.
\end{equation*}
Then, for all large $n$,
\begin{equation}\label{minorationA2}
0  \leq \sum\limits_{i=\c+1}^{m_n-1}\exp(V(i)) \leq 2\cccc
n^{-\e'/2}e^{\d} \leq 1.
\end{equation}
Accordingly (see Zeitouni \cite{Z1} p. 196),
\begin{equation*} \po(\tau_{m_n}<\tau_{-1})
=  \frac{\exp(V(-1))}{\sum_{k=-1}^{m_n-1}\exp(V(k))}
 \geq
\frac{\frac{\e_0}{1-\e_0}}{\frac{1-\e_0}{\e_0}+\frac{1}{1-e^{-\d}}+1}:=C_7>0.
\end{equation*}
\hfill $\Box$

We denote by $\po^x$ and $\eo^x$ the probability and expectation
of $(S_n)_{n\in\N}$, starting at site $x$ and conditioned on the
environment $\o$. We have (see Zeitouni, \cite{Z1}, p. 250)
\begin{fact}
\label{fait} If $a<x<b$,
\begin{equation}
\eo^x(\tau_a    \wedge\tau_b)\leq
\sum_{k=x}^{b-1}\sum_{\ell=a}^{k}\frac{\exp[V(k)-V(\ell)]}{\o_{\ell}}.
\end{equation}
\end{fact}

\bigskip
We can now give an upper bound for the hitting time of $m_n$ if
the RWRE hits $m_n$ before $-1$:
\begin{lemma}\label{lemma34}
There exists a constant $C_8>0$ such that
\begin{equation*}
\forall\o\in \EA,\qquad \po(\tau_{m_n}<\tau_{-1}\textnormal{ and }
\tau_{m_n}\leq n^{4\e'})\geq C_8.
\end{equation*}
\end{lemma}

\bigskip
\noindent {\bf Proof:} According to Fact \ref{fait} and
\eqref{mn}, we obtain for $\o\in \EA$,
\begin{eqnarray*}
\eo(\tb\wedge\tau_{m_n}) & \leq &
\sum_{k=0}^{m_n-1}\sum_{\ell=-1}^{k}
\frac{\exp[V(k)-V(\ell)]}{\o_{\ell}}\\
& \leq & \frac{1}{\e_0} (m_n+1)^2\exp(\cc+2\ddd)\\
& \leq & \frac{1}{\e_0}[(1+\e')(\log n)^2+1]^2 n^{2\e'}\leq
n^{3\e'}.
\end{eqnarray*}
Now, by Chebyshev's inequality,
\begin{eqnarray*}
\po(\tb\wedge\tau_{m_n}\geq n^{4\e'}) \leq
n^{-4\e'}\eo(\tb\wedge\tau_{m_n})\leq n^{-\e'}.
\end{eqnarray*}
As a consequence, recalling Lemma \ref{lemma32},
\begin{eqnarray*}
\po(\tau_{m_n}<\tau_{-1}\text{ and } \tau_{m_n}\leq n^{4\e'})
& = & \po(\tau_{m_n}<\tau_{-1})-\po(n^{4\e'}<\tau_{m_n}<\tb)\\
& \geq & C_7-n^{-\e'}\geq C_7/2:=C_8
\end{eqnarray*}
for $n$ large enough. \hfill$\Box$


\subsection{Leaving the valley}\label{sub32}
First, we give a majoration of the exit time from the valley
$(0,m_n,b_n)$.

\begin{lemma}\label{major} There exists a constant $C_9>0$,
depending only on $\eta$ and $\e$, such that
\begin{equation*}
\forall\o\in \EA,\qquad \po^{m_n}(\tau_{-1}\leq n^{1-5\e'})\geq
C_9.
\end{equation*}
\end{lemma}

\bigskip
\bigskip
\noindent {\bf Proof:} Let $\o\in E_n$. The probability to leave
the valley $(0,m_n,b_n)$ on the left is
\begin{eqnarray}
\po^{m_n}(\tau_{-1}<\tau_{b_n}) & = &
\frac{1}{\frac{\sum_{k=-1}^{m_n-1}\exp(V(k))}
{\sum_{k=m_n}^{b_n-1}\exp(V(k))}+1} \nonumber\\
& \geq &
\frac{1}{\left(1+\frac{1}{1-e^{-\d}}+\frac{1-\e_0}{\e_0}\right)\frac{1-\e_0}{\e_0}+1}
:=2C_9,\label{e1p10}
\end{eqnarray}
due to \eqref{minorationA1} and \eqref{minorationA2}, and since
$\exp(V(b_n))\geq 1$.

Moreover, Fact \ref{fait} gives (by symmetry), recalling
\eqref{mn},
\begin{eqnarray*}
\eo^{m_n}(\tau_{b_n}\wedge\tb) & \leq & \sum_{k=0}^{m_n}\
\sum_{\ell=k}^{b_n}
\frac{\exp[V(k-1)-V(\ell-1)]}{\e_0}\\
& \leq & \frac{1}{\e_0} (b_n+1)^2\frac{1-\e_0}{\e_0}\exp[V(0)-V(m_n)]\\[1mm]
& \leq & (3\log^2 n)^2 \exp[(1-6\e') \log n]\e_0^{-2}\\
& \leq & n^{1-11\e'/2}
\end{eqnarray*}
for $n$ large enough. Then Chebyshev's inequality yields
\begin{equation*}
\po^{m_n}(n^{1-5\e'}<\tau_{b_n}\wedge\tb)\leq n^{-\e'/2}.
\end{equation*}
Consequently, for all environment $\o\in \EA$, recalling
\eqref{e1p10},
\begin{eqnarray*}
\po^{m_n}(\tau_{-1}  \leq  n^{1-5\e'})& \geq &
\po^{m_n}(\tau_{-1}\leq \tau_{b_n} \text{ and } \tau_{-1}\leq
n^{1-5\e'})\\
& \geq & 2C_9-n^{-\e'/2}\geq C_9
\end{eqnarray*}
for $n$ large enough. $\hfill\Box$
\bigskip

Now we give a lower bound for $\tau_{-1}$.

\begin{lemma}\label{minor}
We have,
\begin{equation*}
\inf_{\o\in
E_n}\po^{m_n}(\tau_{-1}>n^{1-10\e'})\underset{n\to\infty}
{\longrightarrow}1.
\end{equation*}
\end{lemma}

\bigskip
\noindent {\bf Proof:} Let $\o\in E_n$. To establish Lemma
\ref{minor}, we use another argument of Sinai's proof.
When the RWRE is located at $m_n-1$, the probability that it hits
$-1$ before going to $m_n$ is
\begin{eqnarray*}
\po^{m_n-1}(\tau_{-1}<\tau_{m_n})
& = & \frac{\exp (V(m_n-1))}{\sum_{k=-1}^{m_n-1}\exp (V(k))}\\
&  \leq & \exp(V(m_n-1)) \\
&  \leq &
\left(\frac{1-\e_0}{\e_0}\right)\frac{e^{\d}}{n^{1-9\e'}},
\end{eqnarray*}
due to \eqref{mn}. Similarly, we have
\begin{equation*}
\po^{m_n+1}(\tau_{b_n}<\tau_{m_n}) \leq
\left(\frac{1-\e_0}{\e_0}\right)\frac{e^{\d}}{n^{1-9\e'}}:=
\frac{C_{10}}{n^{1-9\e'}} .
\end{equation*}
As the RWRE is recurrent, we can consider the $\lfloor
n^{1-10\e'}\rfloor $ first excursions away from $m_n$, which are
independent under $\po$. More precisely, let us define recursively
\begin{equation*}
\left\{
\begin{array}{lcl}
\tau_{m_n}^{(1)} & := & \tau_{m_n},\\[2mm]
\tau_{m_n}^{(k+1)} & := & \inf\{\ell>\tau_{m_n}^{(k)},\quad
S_{\ell}=m_n\}, \quad k\geq 1,
\end{array}
\right.
\end{equation*}
and consider the set
\begin{equation*}
\Gg:=\left\{ \forall 1\leq k\leq \lfloor n^{1-10\e'}\rfloor,\quad
\tau_{m_n}^{(k)}<\tau_{-1}\wedge\tau_{b_n}\right\}.
\end{equation*}
We obtain
\begin{eqnarray*}
\po^{m_n}(\Gg^c)\leq \lfloor
n^{1-10\e'}\rfloor\po^{m_n\pm1}(\tb\wedge\tau_{b_n}<\tau_{m_n})
\leq C_{10}n^{-\e'}.
\end{eqnarray*}
Now, on $\Gg$, the RWRE $(S_i)_{i\geq \tau_{m_n}}$ stays in
$[0,b_n]$ during the first $\lfloor n^{1-10\e'}\rfloor$ excursions
away from $m_n$, hence $\tau_{-1}>n^{1-10\e'}$. Therefore,
\begin{equation*}
\forall\o\in E_n,\qquad \po^{m_n}(\tau_{-1}>n^{1-10\e'})\geq
\po^{m_n}(\Gg)\geq  1-C_{10}n^{-\e'}.
\end{equation*}
 $\hfill\Box$

Combining Lemmas \ref{major} and \ref{minor}, we get
\begin{equation*}
\forall \o\in E_n,\qquad \po^{m_n}(n^{1-10\e'}<\tau_{-1}\leq
n^{1-5\e'})\geq \frac{C_9}{2}
\end{equation*}
for $n$ large enough. Recalling Lemma \ref{lemma34}, this ends the
proof of Lemma \ref{encadrement}.\hfill$\Box$


\mysection{Proofs of Theorem \ref{th1} and Proposition
\ref{moments}} \label{simpl}

In this section, we use the results of the previous sections to
prove Theorem~\ref{th1} and Proposition~\ref{moments}.

\subsection{Proof of Proposition \ref{moments}}
Let $\a\in\R_+^*$ and
\begin{equation*}
M_{\a}:=\sup_{x\in\R_+}(x^{\a}e^{-x})\in(0,+\infty).
\end{equation*}
Then,
\begin{equation}\label{minoreq}
\forall r>0, \qquad \E[\tb^{\alpha} \exp(-r\tb)]\leq M_\alpha
r^{-\alpha}.
\end{equation}

Now we give a lower bound for $\eo[\tb^\alpha \exp(-r \tb)]$. For
any $0<a<1$ and any $\o$,
\begin{equation*}
\eo[\tb^{\a}\exp(-\tb/n)]\geq e^{-1} n^{\a a} \po(n^a\leq\tb\leq
n).\end{equation*}

\noindent Thus, by Lemma \ref{encadrement}, for any $\e>0$, taking
$a=1-10\e'=1-10\d\e$,
\begin{equation*}
\forall \o\in \EA,\qquad
 \eo[\tb^{\alpha} \exp(-\tb/n)]\geq
C_6 e^{-1}n^{\a(1-10\d\e)}.
\end{equation*}
Integrating this inequality on $\EA$, and in view of Lemma
\ref{en}, we get, for all large $n$,
\begin{equation*}
\E[\tb^{\alpha} \exp(-\tb/n)]\geq C_6 C_0
e^{-1}n^{\a(1-10\d\e)-\theta\e}.
\end{equation*}
Since $\e>0$ can be arbitrary small, this, together with
\eqref{minoreq}, yields
\begin{equation*}
\E(\tau_{-1}^\a
e^{-r\tau_{-1}})=\left(\frac{1}{r}\right)^{\a+o(1)},\qquad r\to
0^+.
\end{equation*}
By symmetry, we can replace $\tau_1$ by $\tau_{-1}$, which gives
Proposition \ref{moments}.\hfill$\Box$


\subsection{Proof of Theorem \ref{th1}}
\indent It is known (see den Hollander \cite{H3}, p. 80), that
(\ref{pb1}) is equivalent to
\begin{equation}
\lim_{r\ra 0^-} \frac{[\log \l]''(r)}{\{[\log \l]'(r)\}^3}=0,
\label{pb2}
\end{equation}
where
\begin{equation}
\log\lambda(r)=\E[\log\eo(e^{r\tau_{1}})].
\end{equation}
Note that
\begin{eqnarray*}
f(r) & := &  \frac{[\log \l]''(r)}{\{[\log \l]'(r)\}^3}\\
& \leq & \frac{\E\left(\frac{\eo(\tau_1^2
e^{r\tau_1})}{\eo(e^{r\tau_1})}\right)}
{\left[\E\left(\frac{\eo(\tau_1 e^{r
\tau_1})}{\eo(e^{r\tau_1})}\right)\right]^3}:=g(r).
\end{eqnarray*}
Moreover, due to assumption (\ref{lip}), we have, for all $-1<r<0$
and for all environment $\o$,
\begin{equation*}
\e_0 e^{-1}\leq\o_0e^{r}\leq\eo(\exp(r\tau_1))\leq 1.
\end{equation*}
As a consequence, for $-1<r<0$,
\begin{equation*}
g(r)\leq \frac{e}{\e_0}
\frac{\E[\tau_1^2\exp(r\tau_1)]}{\left\{\E[\tau_1
\exp(r\tau_1)]\right\}^3}:=h(r).
\end{equation*}
Furthermore, $f(r)\geq 0$ (by the Cauchy--Schwarz inequality).
 Now, according to Proposition \ref{moments},

$$h(r) = |r|^{1+o(1)}\underset{r\to0^-}{\longrightarrow 0},$$
which proves \eqref{pb2} and thus Theorem \ref{th1}.\hfill$\Box$


\section{Comparison between rate functions}\label{sectbrox}

In this section we consider the diffusion $X$ in the random
potential $W_{\k}$ and assume $\k>1$. In this case,
$v_{\k}=(\k-1)/4$. We know (see Taleb \cite{Ta}) that the rate
function $\jk$ of quenched large deviations for $X$ can be written
as $\jk(x)=x \ik(1/x)$ for $x> 0$, where
\begin{equation}\label{quotient}
\ik(u)=\sup_{\l\geq 0}(\gk(\l)-\l u),
\end{equation}
and $\gk$ can be expressed in terms of modified Bessel functions
(see \eqref{bessel} below).

Let
\begin{equation*}
\fv(\l):=\sqrt{2\l+{v_{\k}}^2}-{v_{\k}}.
\end{equation*}
We first show that $\gk(\l)<\fv(\l)$ for large $\l$. Then we use a
differential equation satisfied by $\gk$ to prove that this
inequality is true on $\R^*_+$. Finally, we prove Theorem
\ref{comptaux} and Proposition \ref{prop14}.


\subsection{Study in the neighbourhood of  $+\infty$}
According to Taleb (we mention that in Taleb~\cite{Ta}, p. 1178,
the expression $F_{\k}(\l)$ should be
$2(2\l)^{\k/2}K_{\k}[4\sqrt{2\l}]$, see for instance Magnus et
al., \cite{MOS} p. 85; this misprint has no consequence on the
results of \cite{Ta}), we have
\begin{equation}\label{bessel}
\forall \l\geq 0,\qquad
\gk(\l)=\sqrt{2\l}\frac{K_{\k-1}(4\sqrt{2\l})}{K_{\k}(4\sqrt{2\l})}.
\end{equation}
Using the ``series of the Hankel type'' (see Magnus et al.
\cite{MOS}, p. 139), we obtain
\begin{equation}\label{dl}
\gk(\l)=\sqrt{2\l}-\frac{1}{4}\left(\k-\frac{1}{2}\right)+O\left(\frac{1}{\sqrt{\l}}\right)
\qquad \l\to+\infty.
\end{equation}
This yields
\begin{equation*}
\gk(\l)-\fv(\l)\underset{\l\to+\infty}{\longrightarrow}-\frac{1}{8}.
\end{equation*}
Consequently, there exists $B>0$, such that
\begin{equation}
\forall\l\geq B,\qquad  \gk(\l)<\fv(\l). \label{vois}
\end{equation}


\subsection{Using a differential equation}
According to Taleb \cite{Ta}, $\gk$ is a solution of the
differential equation $x y'-2y^2-\k y=-4x$ on $(0,+\infty)$. It is
natural to introduce
\begin{equation}\label{eqdefA}
A(x):=x\fv'(x)-2\fv^2(x)-\k\fv(x)+4x=\frac{-x-v_{\k}^2+v_{\k}
\sqrt{2x+v_{\k}^2}}{\sqrt{2x+v_{\k}^2}}.
\end{equation}
In particular, $A(x)<0$ for all $x>0$.

Let us consider the set
\begin{equation*} E:=\{x>0,\quad  \gk(x)\geq
\fv(x)\}.
\end{equation*}
We prove that $E=\emptyset$. Indeed, let us assume that $E\neq
\emptyset$. According to~\eqref{vois}, $E\cap
[B,+\infty)=\emptyset$. Consequently, $E$ would have a supremum
$x_0\in (0,B]$. By continuity, $\gk(x_0)=\fv(x_0)$. Now,
\eqref{eqdefA} would yield
\begin{eqnarray*}
\fv'(x_0) & = & \frac{1}{x_0}[A(x_0)+2\fv^2(x_0)+\k\fv(x_0)-4x_0]\\
          & = & \frac{1}{x_0}[A(x_0)+2\gk^2(x_0)+\k\gk(x_0)-4x_0] \\
          & = & \frac{A(x_0)}{x_0}+\gk'(x_0)<\gk'(x_0).
\end{eqnarray*}
Consequently, there would exist an $\e>0$ such that
\begin{equation*}
\forall x\in [x_0,x_0+\e],\qquad \fv(x)<\gk(x).
\end{equation*}
Therefore, $[x_0,x_0+\e]\subset E$, which contradicts $x_0=\sup
E$. Hence $E=\emptyset$, which means that

\begin{equation} \forall \l> 0,\qquad  \gk(\l)<\fv(\l). \label{Q3}
\end{equation}


\subsection{Proofs of Theorem \ref{comptaux} and Proposition
\ref{prop14}}

It is easily seen that
\begin{equation*} \forall \l\geq 0,\quad
\inf_{0< u< \frac{1}{v_{\k}}} \left\{ \l u+ \frac{u}{2}
\left(\frac{1}{u}-{v_{\k}}\right)^2\right\}
=\sqrt{2\l+{v_{_k}}^2}-{v_{\k}}= \fv(\l). \label{Q4}
\end{equation*}
Thus \eqref{Q3} yields
\begin{equation}
\forall 0<u<\frac{1}{v_{\k}},\ \forall \l> 0,\qquad  \gk(\l)-\l u<
\frac{u}{2}\left(\frac{1}{u}-v_{\k}\right)^2.\label{Q2}
\end{equation}
Notice that \eqref{Q2} remains true for $\l=0$ since $\gk(0)=0$.
Now, fix $u\in(0,1/v_{\k})$. Recalling \eqref{dl}, it follows that
\begin{equation*}
\gk(\l)-\l u\underset{\l\to+\infty}{\longrightarrow}-\infty.
\end{equation*}
As the function $\l\mapsto [\gk(\l)-\l u]$ is continuous on
$\R_+$, it has a maximum on, say,  $\l_u\in \R_+$.
Hence, by \eqref{Q2},
\begin{equation*}
\sup_{\l\geq 0} (\gk(\l)-\l u )=\gk(\l_u)-\l_u u
<\frac{u}{2}\left(\frac{1}{u}-v_{\k}\right)^2,
\end{equation*}
which can be written as, recalling \eqref{quotient}:
\begin{equation*}
\forall 0<u<\frac{1}{v_{\k}},\qquad
\ik(u)<\frac{u}{2}\left(\frac{1}{u}-v_{\k}\right)^2.
\end{equation*}
This is equivalent to
\begin{equation*}
\forall x>v_{\k},\qquad
\jk(x)=x\ik\left(\frac{1}{x}\right)<J_{v_{\k}}^B(x) =
\frac{1}{2}(x-v_{\k})^2,
\end{equation*}
proving Theorem \ref{comptaux}.

\bigskip
We notice that \eqref{Q3} can be written in terms of modified
Bessel functions, using \eqref{bessel}, which gives Proposition
\ref{prop14}.


\subsection{Remarks}

Recall that the rate function of large deviations of the standard
Brownian motion is $x\mapsto x^2/2$. By the same arguments as in
the case $\k>1$, we obtain for the transient case with zero speed
($0<\k\leq1$),

\begin{prop} (zero speed case),\\
\begin{tabular}{ll}
If $\k\in (0,1/2)$, &then $\forall x>0,\ \jk(x)>x^2/2$;\\
If $\k=1/2$, &then $\forall x>0,\ \jk(x)=x^2/2$;\\
If $\k\in (1/2,1]$, &then $\forall x>0,\ \jk(x)<x^2/2$.
\end{tabular}
\end{prop}

\bigskip
(The case $\k=1/2$ was obtained by Taleb, \cite{Ta}).

 We also notice that Proposition \ref{prop14} together with the formula
$K_{\nu-1}(z)-K_{\nu+1}(z)=-\frac{2\nu}{z}K_{\nu}(z)$ also give a
lower bound for $K_\nu/K_{\nu+1}$:
\begin{equation*}
\forall\nu>0, \forall y>0, \qquad \frac{K_{\nu}(y)}{K_{\nu+1}(y)}>
\frac{1}{y}\left[\frac{y^2}{\sqrt{y^2+(\nu+1)^2}-(\nu+1)}-2(\nu+1)\right].
\end{equation*}

\bigskip
\bigskip
\noindent {\bf \Large Acknowledgements}

 I am grateful to Zhan Shi
for several helpful discussions.




\end{document}